\documentclass[12pt,reqno]{amsart}
\usepackage{amssymb,delarray}
\usepackage{amsfonts,amscd}
\usepackage{epsfig}
\usepackage[all]{xy}

\def\leq{\leqslant}
\def\geq{\geqslant}

\newtheorem{thm}{Theorem}
\newtheorem{lem}
{Lemma}
\newtheorem{prop}
{Proposition}
\newtheorem{claim}
{Claim}

\newtheorem{cor}
{Corollary}
{Remark}
{Question}
\newtheorem{ex-thm}{Theorem-Example}

{\catcode`\@=11
\gdef\n@te#1#2{\leavevmode\vadjust{%
 {\setbox\z@\hbox to\z@{\strut#1}%
  \setbox\z@\hbox{\raise\dp\strutbox\box\z@}\ht\z@=\z@\dp\z@=\z@%
  #2\box\z@}}}
\gdef\leftnote#1{\n@te{\hss#1\quad}{}}
\gdef\rightnote#1{\n@te{\quad\kern-\leftskip#1\hss}{\moveright\hsize}}
\gdef\?{\FN@\qumark}
\gdef\qumark{\ifx\next"\DN@"##1"{\leftnote{\rm##1}}\else
 \DN@{\leftnote{\rm??}}\fi{\rm??}\next@}}

\begin{document}
\baselineskip=13.7pt plus 2pt 

\title[On the monodromy of the inflection points of plane curves] {On the monodromy of the inflection points of plane curves}
\author[Vik.S. Kulikov]{Vik.S. Kulikov}

\address{Steklov Mathematical Institute}
 \email{kulikov@mi.ras.ru}

\dedicatory{} \subjclass{}
\thanks{The research was partially supported by grant of RFBR (no.15-01-02158).}


\begin{abstract}
We prove that the monodromy group of the inflection points of plane curves of degree $d$ is the symmetric group $\mathbb S_{3d(d-2)}$ for $d\geq 4$ and in the case $d=3$ this group is the group of the projective transformations of $\mathbb P^2$ leaving invariant the nine inflection points of the Fermat curve of degree three.
\end{abstract}

\maketitle
\setcounter{tocdepth}{1}


\def\st{{\sf st}}


\setcounter{section}{-1}

\section{Introduction.} \label{introduc}
Let  $F(\overline a,\overline z)=\displaystyle \sum_{k+m+n=d}a_{k,m,n}z_1^kz_2^mz_3^n$ be the homogeneous polynomial of degree $d$ in variables $z_1,z_2,z_3$ and of degree one in variables $a_{k,m,n}$, $k+m+n=d$. Denote by $\mathcal C_d\subset \mathbb P^{K_d}\times \mathbb P^2$, where $K_d=\frac{d(d+3)}{2}$, the complete family of plane curves of degree $d$ given by equation $F(\overline a,\overline z)=0$.
Let $\widetilde h_d:\mathcal C_d\to  \mathbb P^{K_d}$ and $h_d:{\mathcal I_d}=\mathcal C_d\cap\mathcal H_d\to \mathbb P^{K_d}$ be the restrictions of the projection $\text{pr}_1: \mathbb P^{K_d}\times \mathbb P^2\to \mathbb P^{K_d}$ to $\mathcal C_d$ and ${\mathcal I_d}$ respectively, where
$$\displaystyle \mathcal H_d=\{ (\overline a,\overline z)\in \mathbb P^{K_d}\times \mathbb P^2\mid\det (\frac{\partial^2 F(\overline a,\overline z)}{\partial z_i\partial z_j})=0\} .$$

It is well-known (see, for example, \cite{BK}) that for a generic point $\overline a_0\in \mathbb P^{K_d}$ the intersection of the curve $C_{\overline a_0}=\widetilde h_d^{-1}(\overline a)$ and its Hessian curve $H_{C_{\overline a_0}}$ given by $\frac{\partial^2 F(\overline a_0,\overline z)}{\partial z_i\partial z_j})=0$ is the set of the inflection points of $C_{\overline a_0}$ containing $3d(d-2)$ points. Therefore for $d\geq 3$ we have $\deg h_d=3d(d-2)$.

Let $\mathcal S_d$ be the subvariety of $\mathbb P^{K_d}$ consisting of the points $\overline a$ such that the curves $C_{\overline a}$ are singular and let $\mathcal M_d$ be the closure of subvariety of $\mathbb P^{K_d}$ consisting of the points $\overline a$ such that for $\overline a\in \mathcal M_d$ the curve $C_{\overline a}$ has a $r$-tuple inflection point with $r\geq 2$, i.e., $C_{\overline a}$ has a smooth point $p$ such that the tangent line $L$ to $C_{\overline a}$ at $p$ and  $C_{\overline a}$ have at $p$ the intersection number $(L,C_{\overline a})_p=r+2\geq 4$. Let $\mathcal B_d=\mathcal S_d\cup \mathcal M_d$ (if $d=3$ then $\mathcal M_3=\emptyset$). Then $h_d: {\mathcal I_d}\setminus h_d^{-1}(\mathcal B_d) \to \mathbb P^{K_d}\setminus \mathcal B_d$ is an unramified covering and therefore it defines a homomorphism
$h_{d*}: \pi_1(\mathbb P^{K_d}\setminus \mathcal B_d, \overline a_0) \to \mathbb S_{3d(d-2)}$ (here $\mathbb S_{3d(d-2)}$ is the symmetric group acting on the set $I_{\overline a}=C_{\overline a_0}\cap \mathcal I_d$). The group $\mathcal G_d=\text{Im}\, h_{d*}$ is called the {\it monodromy group of the inflection points of plane curves of degree $d$.}

The main result of the article is the following
\begin{thm} \label{main} The group $\mathcal G_d=\mathbb S_{3d(d-2)}$ if $d\geq 4$ and
$\mathcal G_3$ is a group of order $216$ isomorphic to the group of the projective transformations of $\mathbb P^2$ leaving invariant the nine inflection points of the Fermat curve of degree three.
\end{thm}

The proof of Theorem \ref{main} is given in Section \ref{mains}. To prove Theorem \ref{main}, some properties of the variety ${\mathcal I}_d$ near $r$-tuple inflection points of curves are investigated in this section. In Section \ref{sing}, we investigate properties of ${\mathcal I}_d$ near  a node of a nodal curve of degree $d$ which will be useful in the further investigations of the variety of the inflection points of plane curves of degree $d$.

\section{Proof of Theorem \ref{main}}\label{mains}
\subsection{On the monodromy of dominant morphisms.} \label{monodr}
\label{loc}
Let $\mathcal B\subset \mathbb P^K$ be a reduced hypersurface in the projective space $\mathbb P^K$. It is well-known that the fundamental group $\pi _1(\mathbb P^K\setminus \mathcal B,p)$ is generated by, so called,
{\it bypasses}  $\gamma_{q,L}$ around $\mathcal B$, that is, elements presented by loops $\Gamma_{q,L}$ of the following form. Let $L \subset \mathbb P^K$ be a germ of a smooth curve intersecting the curve $\mathcal B$ at a point $q\in \mathcal B$, $L\not\subset \mathcal B$, and let $S^1 \subset L$ be a circle of small radius with center at $q$. The right orientation on $\mathbb P^K$, defined by complex structure, defines an orientation on $S^1$ and then $\Gamma_{q,L}$ is a loop consisting of a path $l$ lying in $\mathbb P^K\setminus \mathcal B$ and connecting the point $p$ with some point $q_1\in S^1$, the loop $S^1$ (with right orientation) starting and ending at $q_1$, and return to the point $p$ along the path $l$.

Let $f:X\to \mathbb P^K$ be a dominant morphism, where $X$ is a reduced variety, $\dim X=K$. Then there is a hypersurface $\mathcal B\subset \mathbb P^K$ (called the {\it discriminant locus} of $f$) such that $f: Y=X\setminus f^{-1}(\mathcal B)\to \mathbb P^K\setminus \mathcal B$ is an unramified finite covering. Note that $Y$ is a smooth variety.

Let the degree  of the covering $f:Y\to \mathbb P^K\setminus \mathcal B$ is equal to $n$. Then the covering $f$ defines a homomorphism $f_*:\pi_1(\mathbb P^K\setminus \mathcal B,p)\to \mathbb S_n$ (called the
{\it monodromy} of the covering $f$), where the image $G_f:=f_*(\pi_1(\mathbb P^K\setminus \mathcal B,p))$ (called the {\it monodromy group} of $f$) is a subgroup of the symmetric group $\mathbb S_n$ and it acts on the fibre $f^{-1}(p)=\{ p_1,\dots,p_n\}$ as follows. A loop $\Gamma\subset \mathbb P^K\setminus\mathcal B$ representing an element $\gamma \in\pi_1(\mathbb P^K\setminus \mathcal B,p)$ can be lifted to $Y$ and this lift consists of $n$ paths $\Gamma_1,\dots,\Gamma_n\subset Y$  starting and ending at the points of $f^{-1}(p)$. Therefore this lift defines an action $f_*(\gamma)$ on $f^{-1}(p)$ which sends the start point $p_i\in \Gamma_i$ to the end point of  $\Gamma_i$ for each $i=1,\dots,n$.

The following Lemma is obvious.
\begin{lem} \label{pi} In notations used above, let
\begin{itemize}
\item[$(1)$] $\nu :\widetilde L \to f^{-1}(L)$ be the normalization of
the curve $f^{-1}(L)$, where $L\subset\mathbb P^N$;
\item[$(2)$] the preimage $\widetilde f^{-1}(q)$ consist of $k$ points $q_1,\dots, q_k$, where $\widetilde f=f\circ \nu$ and $q\in L\cap \mathcal B$;
\item[$(3)$] $n_i$ be the ramification index of $\widetilde f$ at $q_i$, $n_1+\dots +n_k=n$.
\end{itemize}
Then $f_*(\gamma_{q,L})=c_1\cdot ...\cdot c_k\in\mathbb S_n$ is the product of $k$ pairwise disjoints cycles $c_i$ of length $n_i$.
\end{lem}

Let $V\subset \mathbb P^K$ be a small neighbourhood (in complex-analytic topology) of a point $q\in \mathcal B$ bi-holomorphic to a polidisk $$\Delta^K=\{ (z_1,\dots, z_K)\in \mathbb C^K \mid |z_j|<\varepsilon\ll 1 \,\text{for}\, j=1,\dots,K\}$$ with center at $q$. The imbedding $i:V\setminus \mathcal B \hookrightarrow \mathbb P^K\setminus \mathcal B$ induces a homomorphism $i_*:\pi_1(V\setminus \mathcal B) \rightarrow \pi_1(\mathbb P^N\setminus \mathcal B)$ and a homomorphism $f_{*,loc}=f_*\circ i_*:\pi_1(V\setminus \mathcal B)\to G_f$ defined by $i_*$ and $f_*$ uniquely up to conjugation in $G$. If a neighbourhood $V$ is small enough then the image
$G_{f, q}:=f_{*,loc}(\pi_1(V\setminus \mathcal B))$ does not depend of $V$ and it is called the {\it local monodromy group} of $f$ at the point $q$.
The following Claim is well-known.
\begin{claim}\label{locgroup} In notations used above, let $q$ be a smooth point of $\mathcal B$ and a curve $L$ intersects with $\mathcal B$ transversally at $q$. Then the local monodromy group $G_{f,q}$ is cyclic and it is generated by $f_*(\gamma_{q,L})$.
\end{claim}

\begin{lem} \label{non-si} Let $Z\subset \Delta^{K+1}$ be a germ of a reduced complex-analytic variety, $\dim Z=K$, in the polidisk
$\Delta^{K+1}=\{ (z_1,\dots, z_{K+1}\in \mathbb C^{K+1} \mid |z_j|<\varepsilon\ll 1 \,\text{for}\, j=1,\dots,n+1\}$, $o=(0,\dots, 0)\in Z$, and let the restriction $f: Z\to \Delta^K$ of the projection $\text{pr}: \Delta^{K+1}\to\Delta^K$,
$\text{pr}: (z_1,\dots, z_{K+1})\mapsto (z_1,\dots, z_K)$, has the following properties:
\begin{itemize}
\item[$(i)$] $f$ is a proper finite holomorphic map, $\deg f=n$;
\item[$(ii)$] the discriminant locus $\mathcal B\subset \Delta^K$ is a smooth hypersurface;
\item[$(iii)$] the local degree $\deg_{o} f=n$;
\item[$(iv)$] there is a germ $L\subset \Delta^K$, $\dim L=1$, such that $L$ meets $\mathcal B$ at the point $o'=f(o)$, $L\not\subset \mathcal B$, and $E=f^{-1}(L)$ is a smooth curve.
\end{itemize}
Then $o$ is a non-singular point of $Z$.
\end{lem}
\proof Without loss of generality, we can assume that $\mathcal B$ is given by $z_1=0$ and,
by Weiershtrass preparation theorem, $Z$ is given by equation of the form
\begin{equation}\label{w1} z_{K+1}^n +\sum _{j=0}^{n-1}\alpha_j(z_1,\dots,z_K)z_{K+1}^j=0,\end{equation}
where $\alpha_j(z_1,\dots,z_K)\in \mathbb C[[z_1,\dots,z_K]]$. By Claim \ref{locgroup}, $\alpha_j(0,z_2,\dots,z_K)=0$ for each
$j=0,\dots, n-1$. Therefore $z_1$ is a divisor of each power series $\alpha_j(z_1,\dots,z_K)$,  $\alpha_j(z_1,\dots,z_K)=z_1^{k_j}\beta_j(z_1,\dots,z_K)$, where $\beta_j(z_1,\dots,z_K)$ is a power series coprime with $z_1$ and $k_j$ is a positive integer.  Let the germ $L$ be given by parametrization
\begin{equation} \label{w2} z_j=\sum_{l=1}^{\infty} c_{j,l}t^l, \quad j=1,\dots, K.\end{equation}
Then the curve $E=f^{-1}(L)$ is given by (\ref{w1}) and (\ref{w2}). Therefore $\beta_0(0,\dots,0)\neq 0$, $c_{1,1}\neq 0$, and $k_0=1$, since
$E$ is a smooth curve at $o$. Now, the smoothness of $Z$ follows from inequality $\beta_0(0,\dots,0)\neq 0$. \qed

\subsection{On $r$-tuple inflection points.} \label{beg} In that follows we shall use the following well-known properties of plane curves of degree
$d\geq 3$ (see, for example \cite{BK}). First of all remind that the Hessian curve $H_{C}$ of a plane curve $C$ is independent on the choice of coordinates; $H_C$ intersects $C$ at the singular and inflection points of $C$ if $C$ does not contain a line as its irreducible component. If a line $L$ is a component of $C$, then $L$ also is a component of $H_C$. Moreover, if $\overline z_0$ is a $r$-tuple inflection point of $C$, then (Theorem 1 in subsection 7.3 in  \cite{BK}\footnote{ In \cite{BK}, Theorem 1 is proved under the additional assumption that there is not a line among the irreducible components of $C$. But, it is easy to see that this assumption is not used in the proof of this Theorem.}) the intersection number $(C,H_C)_{\overline z_0}$ at the point $\overline z_0$ is equal to $r$. Therefore we have

\begin{claim} \label{de} Let $C_{\overline a_0}=C\cup (\cup_{j=1}^k L_j)$ be the union of a curve $C$ of degree $\deg C\geq 3$ and $k$ lines $L_j$ (may be, $k=0$). Let $\{\overline z_1,\dots, \overline z_n \}$ be the set of the inflection points of $C_{\overline a_0}$ which do not lie in $\cup_{j=1}^kL_j$. Then there is  a small (in complex analytic topology) neighbourhood $\mathcal U\subset \mathbb P^{K_d}$ of the point $\overline a_0$ such that $h_d^{-1}(U)$ is the disjoint union of $n+1$ open sets $V_l$, $l=1,\dots,n+1$, such that $(\overline a_0,\overline z_l)\in V_l$ for $l\leq n$. The local degree $\deg_{(\overline a_0,\overline z_l)} h_{d}$ of the covering $h_d$ at a point
$(\overline a_0,\overline z_l)$ is equal to $r_l$ if $\overline z_l$ is a $r_l$-tuple inflection point of $C_{\overline a_0}$.  In particular, if $\overline z_l$ is a simple inflection point (that is, $r_l=1$) and  $U$ is chosen small enough, then  $h_{d|V_l}:V_l\to U=h_d(V_l)$ is a bi-holomorphic map.
\end{claim}

\proof The local degree $\deg_{(\overline a_0,\overline z_0)} h_{d}$ of the covering $h_d$ at the point $(\overline a_0,\overline z_0)$ is equal to the intersection number $({\mathcal I},\text{pr}_1^{-1}(\overline a_0))_{(\overline a_0, \overline z_0)}$ of the variety ${\mathcal I}$ and the fibre $\text{pr}_1^{-1}(\overline a_0)$ at $(\overline a_0,\overline z_0)$. On the other hand,
$({\mathcal I},\text{pr}_1^{-1}(\overline a_0))_{(\overline a_0, \overline z_0)}$ is equal to the intersection number of $C_{\overline a_0}$ and its Hessian $H_{C_{\overline a_0}}$ at $\overline z_0$ in $\mathbb P^2$. \qed \\

The following Proposition is well-known, but since I do not know a good reference, a proof will be given.
\begin{prop} \label{prop1} The variety $\mathcal M_d$ is an irreducible hypersurface in $\mathbb P^{K_d}$ for each $d\geq 4$. There is a non-empty Zariski open neighbourhood $\mathfrak M_d\subset  \mathcal M_d$ such that for each $\overline a\in \mathfrak M_d$ the curve
$C_{\overline a}$ is non-singular and it has $3d(d-2)-1$ inflection points {\rm (}that is, it has the only one multiple $(r=2)$ inflection point{\rm )}.
\end{prop}

\proof
Denote by $\mathcal D_r\subset \mathcal C_d$, $r=2,3$, the subfamilies of curves given by
\begin{equation} \label {r=2}
z_1S(z-1,z_2,z_3)+z_2^{r+2}R(z_2,z_3)=0,
\end{equation}
where $S(z_1,z_2,z_3)$ is the generic homogeneous polynomial of degree $d-1$ in variables $z_1,z_2,z_3$ and $R(z_2,z_3)$ is the  generic homogeneous polynomial of degree $d-(r+2)$ in variables $z_2,z_3$. Denote also by $D_r=\text{pr}_1(\mathcal D_r)\subset \mathbb P^{K_d}$ the image of $\mathcal D_r$ under the projection $\text{pr}_1$.
Obviously, $D_{3}\subset D_2$, $D_2$ and  $D_3$ are irreducible projective  varieties, and it is easy to see
that $\dim D_2=\frac{(d-1)(d+2)}{2}+(d-3)=K_d-4$ and $\dim D_3=\frac{(d-1)(d+2)}{2}+(d-4)=K_d-5$.

Similarly, let $\mathcal D_{2,2}\subset \mathcal D_2$ be the subfamily of curves given by
\begin{equation} \label {r=2,2}
z_1z_3S(z_1,z_2,z_3)+z_2^{4}(R_1(z_1,z_2)+R_2(z_2,z_3))=0,
\end{equation}
where $\deg S(z_1,z_2,z_3)=d-2$  and $\deg R_1(z_1,z_2)=\deg R_2(z_2,z_3)=d-4$;  $\mathcal D_{2,2,1}\subset \mathcal D_2$ the subfamily of curves given by
\begin{equation} \label {r=2,2,1}
z_1z_3S(z_1,z_2,z_3)+z_2^{4}z_3R_1(z_2,z_3)+z_1^4R_2(z_1,z_2))=0,
\end{equation}
where $\deg S(z_1,z_2,z_3)=d-2$, $\deg R_1(z_1,z_2)=d-5$, and $\deg R_2(z_2,z_3)=d-4$; and
$\mathcal D_{2,2,2}\subset \mathcal D_2$ the subfamily of curves given  by
\begin{equation} \label {r=2,2,2}
z_1S(z_1,z_2,z_3)+z_2^{4}z_3^4R_1(z_2,z_3)=0,
\end{equation}
where $\deg S(z_1,z_2,z_3)=d-1$ and $\deg R_1(z_2,z_3)=d-8$.
As above, denote  by
$D_{2,2}=\text{pr}_1(\mathcal D_{2,2})$, $D_{2,2,1}=\text{pr}_1(\mathcal D_{2,2,1})$, and $D_{2,2,2}=\text{pr}_1(\mathcal D_{2,2,2})$ the images respectively of $\mathcal D_{2,2}$, $\mathcal D_{2,2,1}$, and $\mathcal D_{2,2,2}$ under the projection $\text{pr}_1$.
Obviously, $D_{2,2}$, $D_{2,2,1}$, and  $D_{2,2,2}$ are irreducible projective  varieties, and it is easy to see that $\dim D_{2,2}=K_d-7$ and $\dim D_{2,2,1}=\dim D_{2,2,2}=K_d-8$.

Let $p_1$ be a $r$-tuple point, $r=2$ or $r\geq3$, of a curve $C$. We choose homogeneous coordinates $(z_1,z_2,z_3)$ so that the point $p_1$ has coordinates $(0,0,1)$ and the line $L_1=\{ z_1=0\}$ is the tangent line to $C$ at the point $p_1$. Then an equation of $C$ has  the form (\ref{r=2}).

Let $C$ has two $2$-tuple inflection points $p_1$ and $p_2$, and let $L_1$ and $L_2$ be the tangent lines to $C$ respectively at $p_1$ and $p_2$. We have three possibilities depending on the position of the points $p_1$ and $p_2$ with respect to the curve $C$:
either $p_1\not\in L_2$ and $p_2\not\in L_1$, or $p_2\in L_1$, but $L_1\neq L_2$, or $L_1=L_2$. In the first case, if we choose homogeneous coordinates $(z_1,z_2,z_3)$ so that the point $p_1$ has coordinates
$(0,0,1)$ and the line $L_1=\{ z_1=0\}$ is the tangent line to $C$ at the point $p_1$, the point $p_2$ has coordinates
$(1,0,0)$ and the line $L_2=\{ z_3=0\}$ is the tangent line to $C$ at the point $p_2$, then an equation of $C$ has  the form (\ref{r=2,2}). In the second case, if we choose homogeneous coordinates $(z_1,z_2,z_3)$ so that the point $p_1$ has coordinates
$(0,0,1)$ and the line $L_1=\{ z_1=0\}$ is the tangent line to $C$ at the point $p_1$, the point $p_2$ has coordinates
$(0,1,0)$ and the line $L_2=\{ z_3=0\}$ is the tangent line to $C$ at the point $p_2$, then an equation of $C$ has  the form (\ref{r=2,2,1}). In the third  case, if we choose homogeneous coordinates $(z_1,z_2,z_3)$ so that the point $p_1$ has coordinates
$(0,0,1)$, the point $p_2$ has coordinates $(0,1,0)$, and the line $L_1=\{ z_1=0\}$ is the tangent line to $C$ at the points $p_1$ and $p_2$, then an equation of $C$ has  the form (\ref{r=2,2,2}).

The group $PGL(3,\mathbb C)$ acts on $\mathbb P^{K_d}$ so that $g(\overline a)\in \mathcal M_d$ for each $g\in PGL(3,\mathbb C)$ and $\overline a\in D_2$. Denote by $\Gamma_{\dag}\subset PGL(3,\mathbb C)$ the subgroup leaving invariant the variety $D_{\dag}$, where
$\dag $ is either $2$, or $3$, or $\{ 2,2\}$, or $\{ 2,2,1\}$, or $\{ 2,2,2\}$.

Obviously, the groups $\Gamma_r$, $r=2$ or $r\geq 3$, contain a subgroup $\Gamma_0$  consisting  of the elements of the following form
\begin{equation}\label{grou}
g=\left( \begin{matrix} g_{1,1}, & 0, & 0 \\
g_{2,1}, & g_{2,2}, & 0 \\
g_{3,1}, & g_{3,2}, & g_{3,3} \end{matrix} \right)
\end{equation}
and it is easy to see that $\Gamma_0$ is a subgroup of finite index in  $\Gamma_r$, since each non-singular curve $C$ can have only finitely many multiple inflection points. Therefore $\dim \Gamma_r=5$. Note that if there exists a curve $C_{\overline a}$, $\overline a\in \mathcal M_d$, which has only one multiple inflection point, then $\Gamma_0=\Gamma_2$.

Similarly, the diagonal group $\Delta$ is a subgroup of $\Gamma_{2,2}$ of finite index; the group consisting of the elements of the form
\begin{equation}\label{grou2}
g=\left( \begin{matrix} g_{1,1}, & 0, & 0 \\
g_{2,1}, & g_{2,2}, & 0 \\
0, & 0, & g_{3,3} \end{matrix} \right) \in PGL(3,\mathbb C)
\end{equation}
is a subgroup of $\Gamma_{2,2,1}$ of finite index; and the group consisting of the elements of the form
\begin{equation}\label{grou3}
g=\left( \begin{matrix} g_{1,1}, & 0, & 0 \\
g_{2,1}, & g_{2,2}, & 0 \\
g_{3,1}, & 0, & g_{3,3} \end{matrix} \right) \in PGL(3,\mathbb C)
\end{equation}
is also a subgroup of $\Gamma_{2,2,1}$ of finite index. Therefore $\dim \Gamma_{2,2}=2$, $\dim \Gamma_{2,2,1}=3$, and $\dim \Gamma_{2,2,2}=4$.

Consider the morphism $\nu: PGL(3,\mathbb C)\times D_2\to \mathbb P^{K_d}$ given by $\nu((g,\overline a))=g(\overline a)$. Obviously, the image $\text{Im}\, \nu\subset \mathcal M_d$ is an everywhere dense subset of $\mathcal M_d$. The preimage of a point $\nu((g_0,\overline a_0))$, where $\overline a_0$ is a generic point of  $D_2$, is the variety $\{ (g, g^{-1}g_0(\overline a_0) \mid g\in \Gamma_r \}$ of dimension five. Therefore $$\dim \nu(PGL(3,\mathbb C)\times D_r)=\dim PGL(3,\mathbb C)+\dim D_r-\dim \Gamma_r=K_d-r+1.$$
In particular, $\dim \mathcal M_d=\dim \nu(D_2)=K_d-1$ and therefore $\mathcal M_d$ is a hypersurface in $\mathbb P^{K_d}$.

Similar calculations give rise $\dim \nu(PGL(3,\mathbb C)\times D_{2,2})=K_d-2$,
$$\dim \nu(PGL(3,\mathbb C)\times D_{2,2,1})=K_d-3,\quad
\dim \nu(PGL(3,\mathbb C)\times D_{2,2,2})=K_d-4.$$
It is well-known that $\mathcal S_d$ is a divisor in $\mathbb P^{K_d}$
and it is easy to see that $\mathcal M_d\not\subset \mathcal S_d$.
Therefore $\dim \mathcal M_d\cap \mathcal S_d=K_d-2$. Now,
$$\mathfrak M_d=\mathcal M_d\setminus (\overline{\nu (PGL(3,\mathbb C)\times (D_{2,2}\cup D_{2,2,1}\cup D_{2,2,2})}\cup \mathcal S_d)$$
is the desired variety, where by $\overline M$ is denoted the closure of a variety $M\subset \mathbb P^{K_d}$.

Therefore $\Gamma_2=\Gamma_0$ and hence $\nu(PGL(3,\mathbb C)\times D_{2})$ is irreducible, since
there is a point $\overline a\in D_2$ such that $C_{\overline a}$ is non-singular and it has only one multiple inflection point (more precisely, $2$-tuple inflection point).
\qed
\subsection{On a quasi-imbedding of the permutation group $\mathcal G_d$ into $\mathcal G_{d+1}$.} \label{quasi}
Denote by $(G,n)$ a subgroup $G$ of the symmetric group $\mathbb S_n$ acting on a set $J_n$ of cardinality $n$ as permutations and call $(G,n)$ a {\it permutation group}.

Let $J_n= J_{m}\bigsqcup J_{k}$, $m+k=n$, be a partition of $J_n$ and $(G,n)$ a permutation group such that the action of $G$ on $J_n$ leaves invariant the set $J_{m}$. The action of $G$ on $J_{m}$ defines a homomorphism
$\varphi_{n,m}: G\to \mathbb S_{m}$, that is, it defines the permutation group $(G_{J_m}, m)$, where $G_{J_m}=Im\, \varphi_{n,m}$. We say that a permutation group $(H,m)$ is {\it quasi-imbedded} in a permutation group $(G,n)$ (and denote this quasi-imbedding by $(H,m)\prec (G,n)$) if
\begin{itemize}
\item[$(i)$] $n\geq m$ and there is a partition $J_n=J_m\bigsqcup J_{n-m}$ such that
$G$ leaves invariant the set $J_m$;
\item[$(ii)$]  the permutation groups $(G_{J_m},m)$ and $(H,m)$ are isomorphic as permutation groups.
\end{itemize}

Remind that the group $\mathcal G_d$ is the image of $\pi_1(\mathbb P^{K_d}\setminus \mathcal B_d,\overline a_0)$ under the homomorphism
$h_{d*}: \pi_1(\mathbb P^{K_d}\setminus \mathcal B_d,\overline a_0)\to\mathbb S_{3d(d-2)}$, where the symmetric group
$\mathbb S_{3d(d-2)}$ acts on $J_{3d(d-2)}:=h_d^{-1}(\overline a_0)$. Therefore in what follows, the group $\mathcal G_d$ will be considered as a permutation group $(\mathcal G_d,3d(d-2))$, but it will be denoted again simply by $\mathcal G_d$.

\begin{claim}\label{imbedding} For each $d\geq 3$  there is a quasi-imbedding of $\mathcal G_d$ in $\mathcal G_{d+1}$.
\end{claim}
\proof For each $g\in \mathcal G_d$ let us choose a continuous loop $\Gamma: [0,1]\to \mathbb P^{K_d}\setminus \mathcal B_d$
representing an element $\gamma \in h_{d*}^{-1}(g)\subset\pi_1(\mathbb P^{K_d}\setminus \mathcal B_d,\overline a_0)$, where
$[a,b]=\{ t\in \mathbb R \mid a\leq t\leq b\}$ is a segment in $\mathbb R$. Then the action of $g$ on $J_{3d(d-2)}$ is defined by the disjoint union $h_{d}^{-1}(\Gamma([0,1]))=\bigsqcup_{j=1}^{3d(d-2)}l_j([0,1])$ of $3d(d-2)$ continuous paths $l_j:[0,1]\to \mathbb P^{K_{d+1}}\setminus \mathcal B_{d+1}$ starting and ending at the points of $h_{d}^{-1}(\overline a_0)$. Denote $\overline z_{\tau,j}=\text{pr}_2(l_j(\tau)$.

Since $h_{d}^{-1}(\Gamma([0,1]))$ is one-dimensional as a topological space, we can choose a line
$L\subset \mathbb P^2$ such that $L\cap\text{pr}_2(\bigsqcup_{j=1}^{3d(d-2)}l_j([0,1]))=\emptyset$. For each $\overline a\in \mathbb P^{K_d}$ the curve $C_{\overline a}\cup L$ has degree $d+1$. Therefore the choice of $L$ defines an imbedding
$\lambda_d :\mathbb P^{K_d}\hookrightarrow \mathbb P^{K_{d+1}}$ given by
$\lambda_d:\overline a \mapsto \overline b =\widetilde h_{d+1}(C_{\overline a}\cup L)\in \mathbb P^{K_{d+1}}$
for $\overline a\in \mathbb P^{K_d}$. 

Consider the loop $\lambda_d(\Gamma([0,1]))$. By Claim \ref{de}, for each $\tau\in [0,1]$ there is  a small (in complex analytic topology) connected neighbourhood $\mathcal U_{\tau}\subset \mathbb P^{K_{d+1}}$ of the point $\lambda_d(\Gamma(\tau))$ such that $h_{d+1}^{-1}(U_{\tau})$ is the disjoint union of $3d(d-2)+1$ open sets $V_{\tau,j}$, $j=1,\dots,3d(d-2)+1$, such that $(\lambda_d(\Gamma(\tau)),\overline z_{\tau,j})\in V_{\tau,j}$ and $h_{d+1}: V_{\tau,j}\to U_{\tau}$ is a bi-holomorphic map for $j\leq 3d(d-2)$ , where $\{ \overline z_{\tau,1},\dots,\overline z_{\tau, 3d(d-2)}\}$ is the set of the inflection points of the curve $C_{\Gamma(\tau)}$. Note that we can choose $U_{0}$ and $U_1$ such that $U_0=U_1$ and this open set does not depend on $g\in \mathcal G_d$.

Let $\Delta_{\tau}=\{ t\in [0,1]\mid \tau-\varepsilon_{\tau}< t<\tau+\varepsilon_{\tau}\}$
be  segments in $[0,1]$ such that $\lambda_d(\Gamma(\Delta_{\tau}))\subset U_{\tau}$ for $0<\tau <1$ and, similarly, $\Delta_0=\{ t\in [0,1]\mid t<\varepsilon_0\}$ and $\Delta_1=\{ t\in [0,1]\mid t>1-\varepsilon_1\}$ be segments such that $\lambda_d(\Gamma(\Delta_{0}))\subset U_{0}$ and $\lambda_d(\Gamma(\Delta_{1}))\subset U_{1}$.

Consider a path $\Theta: [0,1] \to \mathbb P^{K_{d+1}}\times [0,1]$ given by $\Theta : \tau\mapsto (\lambda(\Gamma(\tau)),\tau)$. Obviously, $\{ U_{\tau}\times \Delta_{\tau} \}_{\tau\in [0,1]}$ is a cover of the path $\Theta([0,1])$. Since $\Theta([0,1])$ is a compact, we can choose a finite cover $$\{ U_{0}\times \Delta_{0}, U_{\tau_1}\times \Delta_{\tau_1}, \dots, U_{\tau_n}\times \Delta_{\tau_n},
U_{1}\times \Delta_{1}\}, \quad 0=\tau_0<\tau_1<\dots <\tau_n< \tau_{n+1}=1.$$

It is easy to see that $U_{\tau_j}\cap U_{\tau_{j+1}}\neq \emptyset$ for each $j=0,\dots, n$. Let us choose a point $\overline b_0\in U_0\setminus \mathcal B_{d+1}$ and points $\overline b_{j,j+1}\in (U_{\tau_j}\cap U_{\tau_{j+1}})\setminus \mathcal B_{d+1}$ for
$j=0,\dots, n$. Each variety $U_{\tau_j}\setminus \mathcal B_{d+1}$ is connected. Therefore for $0\leq j\leq n+1$ we can connect the point
$\overline b_{j-1,j}$ with $\overline b_{j,j+1}$ by a continuous path $\Gamma_{j}\subset U_{\tau_j}\setminus \mathcal B_{d+1}$, where
$\overline b_{-1,0}=\overline b_{n+1,n+2}=\overline b_0$.

Consider the set $J_{3(d^2-1)}=h_{d+1}^{-1}(\overline b_0)=\{ \widetilde q_1,\dots,\widetilde q_{3d(d-1)}\dots ,\widetilde q_{3(d^2-1)}\}$, where
$\widetilde q_j\in V_{0,j}$ for $j=1,\dots, 3d(d-1)$. Denote $\widetilde J_{3d(d-2)}=\{ \widetilde q_1,\dots, \widetilde q_{3d(d-2)}\}\subset J_{3(d^2-1)}$. The consecutive join of the paths $\Gamma_j$, $j=0,\dots,n+1$, is a continuous loop
$$\widetilde{\Gamma}=\Gamma_0\cup\dots\cup\Gamma_{n+1}\subset \mathbb P^{K_{d+1}}\setminus \mathcal B_{d+1}$$
starting and ending at $\overline b_0$. Then the start and end points of the paths $\widetilde l_j$, entering in the disjoint union
$h_{d+1}^{-1}(\widetilde{\Gamma})=\bigsqcup_{j=1}^{3(d^2-1)}\widetilde l_j$ of $3(d^2-1)$ continuous paths, are contained in $J_{3(d^2-1)}$. Let $\widetilde g=h_{d+1*}(\widetilde{\gamma})\in \mathcal G_{d+1}$, where $\widetilde{\gamma}\in \pi_1(\mathbb P^{K_{d+1}}\setminus \mathcal B_{d+1},\overline b_0)$ is represented by $\widetilde{\Gamma}$.
If we number the paths $\widetilde l_j$ so that the start point of $\widetilde l_j$ is $\widetilde q_j$ for $j\leq 3d(d-2)$, then it  easily follows from the construction of $\widetilde{\Gamma}$ that
\begin{itemize}
\item[$(1)$] the end point of $\widetilde l_j$ lies also in $\widetilde J_{3d(d-2)}$ for each $j\leq 3d(d-2)$;
\item[$(2)$] $\widetilde g$ leaves invariant the set $\widetilde J_{3d(d-2)}$;
\item[$(3)$] the restriction of the action of $\widetilde g$ to $\widetilde J_{3d(d-2)}$ is the same as the action of $g$ on $J_{3d(d-2)}$ if we identify $\widetilde J_{3d(d-2)}$ with $J_{3d(d-2)}$.
\end{itemize}

Denote by $\widetilde{\mathcal G_d}$ a permutation subgroup of $\mathcal G_{d+1}$ generated by the elements $\widetilde g$, where
$g\in \mathcal G_d$. Obviously, the permutation subgroup $\widetilde {\mathcal G_d}$ defines a quasi-imbedding of $\mathcal G_d$ in
$\mathcal G_{d+1}$. \qed

\subsection{Behaviour of the covering $h_d$ near the inflection points of the Fermat curve.} \label{hF}
Let $F_d\subset\mathbb P^2$ be the Fermat curve of degree $d$, i.e., the curve given by equation
$z_1^d+z_2^d+z_3^d=0$. It has $3d$ the $(d-2)$-tuple inflection points
$\overline z_{j,l}$, where
$$\overline z_{1,l}=(0,\mu_l,1),\,\,  \overline z_{2,l}=(\mu_l,0,1),\,\, \overline z_{3,l}=(\mu_l,1,0),\,\, l=1,\dots, d,\,\,
\mu_l=e^{\pi i(2l-1)/d}.$$

The subgroup $G_d$ of $PGL(3,\mathbb C)$, generated by
$$ g_1=\left( \begin{matrix}
0, & 1, & 0 \\ 1, & 0, & 0 \\
0, & 0, & 1 \end{matrix}\right),  g_2=\left( \begin{matrix}
 1, & 0, & 0 \\   0, & 0, & 1 \\
0, & 1, & 0 \end{matrix}\right),  g_3=\left( \begin{matrix}
e^{2\pi i/d}, & 0, & 0 \\ 0, & 1, & 0 \\
0, & 0, & 1 \end{matrix}\right),
$$
leaves invariant the curve $F_d$ and acts transitively on the set of its inflection points. As a group,
it is isomorphic to $\mathbb Z_d^2\ltimes \mathbb S_3$.

Let $f_d\in \mathbb P^{K_d}$ be the point corresponding to the Fermat curve $F_d$ and let $\mathbb C^{K_d}\subset \mathbb P^{K_d}$ be the affine space given by $a_{0,0,d}\neq 0$. Denote again the non-homogeneous coordinates in $\mathbb C^{K_d}$ by $a_{k,m,n}$, $(k,m,n)\neq (0,0,d)$ (here we assume that $a_{0,0,d}=1$). Then the point $f_d\in \mathbb C^{K_d}$ has the following coordinates: $a_{d,0,0}=a_{0,d,0}=1$, and all other coordinates are equal to zero.

Consider a neighbourhood
$$\begin{array}{ll}U_{\varepsilon}=\{ \overline a\in\mathbb C^{K_d} \mid & |a_{k,m,n}-1|<\varepsilon, \,\, \text{for}\,\,\, (k,m,n)= (d,0,0)\,\, \text{or}\,\, (0,d,0)\,\, \text{and} \\ & |a_{k,m,n}|<\varepsilon \,\,\, \text{for all}\,\,\, (k,m,n)\neq (d,0,0)\,\, \text{or}\,\, (0,d,0)\}\end{array}$$
 of the point $f_d$.

\begin{claim} \label{cl1} For positive $\varepsilon \ll 1$ the variety $\mathcal U_{\varepsilon}=h_d^{-1}(U_{\varepsilon})$ is the disjoint union of $3d$ irreducible varieties $\mathcal U_{j,l}$, where for each $j=1,2,3$ and $l=1,\dots,d$  the variety $\mathcal U_{j,l}$ is defined by the following
condition: the $(d-2)$-tuple inflection point $\overline z_{j,l}$ of the curve $F_d$ lies in $\mathcal U_{j,l}$. The restriction $\mathcal U_{j,l}\to U_{\varepsilon}$
of the morphism $h_d$ to each $\mathcal U_{j,l}$ has degree $d-2$.
\end{claim}
\proof Obviously, Claim \ref{cl1} is true in the case $d=3$. Therefore we will assume that $d\geq 4$.

It is easy to see that if $\varepsilon$ is small enough, then the variety $\mathcal U_{\varepsilon}=h_d^{-1}(U_{\varepsilon})=\bigsqcup_{j=1}^3\bigsqcup_{l=1}^d\mathcal U_{j,l}$ is the disjoint union of $3d$ varieties $\mathcal U_{j,l}$ such that $(f_d,\overline z_{j,l})\in \mathcal U_{j,l}$. Therefore to prove Claim \ref{cl1}, it suffices to prove that $\mathcal U_{1,1}$ is an irreducible variety, since the induced actions of the group $G_d$ on $\mathbb P^{K_d}\times \mathbb P^2$ and $\mathbb P^{K_d}$ leave invariant the varieties $U_{\varepsilon}$, $\widetilde h_d^{-1}(U_{\varepsilon})$, and $\mathcal U_{\varepsilon}$, $g\circ h_d=h_d\circ g$ for each $g\in G_d$, and this action induces a transitive action on the set of varieties $\mathcal U_{j,l}$. Obviously, the restriction of $h_d$ to each $\mathcal U_{j,l}$ has degree $d-2$.

To prove that $\mathcal U_{1,1}$ is an irreducible variety, consider a family $C_u$ of curves in $\mathbb P^{K_d}\times \mathbb P^2$
given by
\begin{equation} \label{curve1} F(u,\overline z):=z_1^d+z_2^d+z_3^d+uz_1^2z_3^{d-2}=0\end{equation}
and its image $L=\widetilde h(C_u)\simeq \mathbb C$. Denote by $L_{\varepsilon}=L \cap U_{\varepsilon}$. The family $C_u$ lies in
$L\times \mathbb P^2\subset \mathbb P_{K_d}\times \mathbb P^2$. It is easy to check that in coordinates $(u;z_1,z_2,z_3)$ the Hessian of the family $C_u$ is
\begin{equation} \label{curve2} \displaystyle \mathcal H(u;\overline z)=\det \left(\frac{\partial^2 F(u,\overline z)}{\partial z_i\partial z_j} \right)= d(d-1)z_2^{d-2}z_3^{d-4}H(u,z_1,z_3),\end{equation}
where
$$\begin{array}{l} H(u,z_1,z_3)  =  \\ (d(d-1)z_1^{d-2}+2uz_3^{d-2})(d(d-1)z_3^2 +(d-2)(d-3)uz_1^2)-4(d-2)^2u^2z_1^2z_3^{d-2}=
\\
\displaystyle \frac{d!}{(d-4)!}u z_1^d+d^2(d-1)^2z_1^{d-2}z_3^2-2(d-1)(d-2)u^2z_1^2z_3^{d-2}+2d(d-1)uz_3^d.
\end{array} $$

Therefore the curve $J=h_d^{-1}(L)\subset \mathbb C\times \mathbb P^2$, given by $F(u,\overline z)=\mathcal H(u,\overline z)=0$, is the union of curves, $J=J_1\cup J_2\cup J_3$ (if $d=4$ then $J_3=\emptyset$), where $J_2$ and $J_3$ (if $d\geq 5$)  are the
intersections of the surface given by (\ref{curve1}) and, respectively, two surfaces given by $z_2=0$ and $z_3=0$, and $J_1$ is the intersection of the surface given by equation (\ref{curve1}) and the surface $\overline H_1$ given by
\begin{equation} \label{curve3} \displaystyle \frac{d!}{(d-4)!}u z_1^d+d^2(d-1)^2z_1^{d-2}z_3^2-2(d-1)(d-2)u^2z_1^2z_3^{d-2}+2d(d-1)uz_3^d=0. \end{equation}

It is easy to see that $(\cup_{l=1}^d\mathcal U_{1,l})\cap (J_2\cup J_3)=\emptyset$ and
$(\cup_{l=1}^d\mathcal U_{1,l})\cap J_1\subset \mathbb C \times \mathbb C^2$, where
$\mathbb C^2=\mathbb P^2\setminus \{z_3=0\}$. Let $x=z_1/z_3$, $y=z_2/z_3$ be coordinates in $\mathbb C^2$, then the surface $H_1= \overline H_1 \cap (\mathbb C \times \mathbb C^2)$  is given by  equation
\begin{equation} \label{curve4} \displaystyle \frac{d!}{(d-4)!}u x^d+d^2(d-1)^2x^{d-2}-2(d-1)(d-2)u^2x^2+2d(d-1)u=0. \end{equation}
Since the polynomial in equation (\ref{curve4}) depends only on the variables $x$ and $u$, the surface $H_1$ is isomorphic to the product $E\times \mathbb C^1$, where $E$ is a curve in $\mathbb C^2$ given by equation (\ref{curve4}).

Let us show that the polynomial $H(u,x)$ in the left side of (\ref{curve4}) is irreducible in the ring $\mathbb C[u,x]$. Indeed, assume that
$H(u,x)=H_1(u,x)H_2(u,x)$. Then $H_1(u,x)=A_1(x)u+A_2(x)$ and $H_2(u,x)=A_3(x)u +A_4(x)$, since $H(u,x)$ is a polynomial of degree two in variable $u$ and the polynomials $2(d-1)(d-2)x^2$ and $\frac{d!}{(d-4)!} x^d+2d(d-1)$ are coprime. Therefore we have
\begin{equation} \label{equ1} A_1(x)A_3(x)=-2(d-1)(d-2)x^2,\quad A_2(x)A_4(x)=d^2(d-1)^2x^{d-2},\end{equation}
\begin{equation}\label{equ2} A_1(x)A_4(x)+A_2(x)A_3(x)=\frac{d!}{(d-4)!} x^d+2d(d-1).\end{equation}
It follows from (\ref{equ1}) that $A_1(x)=b_1x^{k_1}$, and  $A_3(x)=b_3x^{2-k_1}$, where $0\leq k_1\leq 2$ and
$b_1b_3=-2(d-1)(d-2)\in\mathbb C$. Similarly, $A_2(x)=b_2x^{k_2}$ and  $A_4(x)=b_4x^{d-2-k_2}$, where $0\leq k_2\leq d-2$ and $b_2b_4=d^2(d-1)^2\in\mathbb C$. Therefore
\begin{equation}\label{equ3} b_1b_2b_3b_4=-2d^2(d-1)^3(d-2). \end{equation}
It follows from (\ref{equ2}) that
$$b_1b_4x^{d+k_1-k_2-2}+ b_2b_3x^{k_2+2-k_1}=\frac{d!}{(d-4)!} x^d+2d(d-1)$$
and therefore
$$b_1b_4=\frac{d!}{(d-4)!}\,\, \text{and}\,\, b_2b_3=2d(d-1)\,\, \text{or}\,\, b_1b_4=2d(d-1)\,\, \text{and}\, b_2b_3=\frac{d!}{(d-4)!},$$
but in both cases we have
\begin{equation}\label{equ4} b_1b_2b_3b_4=2\frac{d!}{(d-4)!}d(d-1)=2d^2(d-1)^2(d-2)(d-3). \end{equation}
It follows from (\ref{equ3}) and (\ref{equ4}) that we have the equality
$$2d^2(d-1)^2(d-2)(d-3)=2d^2(d-1)^2(d-2)(d-3),$$
i.e., $d=0,1$ or $2$, but, by assumption, $d\geq 3$ and therefore $H(u,x)$ is an irreducible polynomial.

Denote by $S$ the union of the set of critical values of the restriction $p:E\to \mathbb C\simeq L$  of the projection $\text{pr}\, :(u,x)\mapsto (u)$ to the irreducible curve $E$ and the set of the images under the projection of the intersection points of $E$ and the curve given by $x^d+ux^{2}+1=0$. Note that $S$ is a finite set. Let $S=\{ 0, u_1,\dots, u_t\}$. Then for $\varepsilon<< 1$ such that $\varepsilon < \min u_s$, where minimum is taken over all $u_s\in S\setminus \{ 0\}$, and for each fixed non-zero value $u_0$ of $u$, $|u_0|<\varepsilon$, the set $p^{-1}(u_0)=\{ (u_0,x_1(u_0)),\dots, (u_0,x_d(u_0))\}$ consists of $d$ different points such that two of these points, say $(u_0,x_{d-1}(u_0))$ and $(u_0,x_d(u_0))$, lie very far from the point $(0,0)$ and the other $d-2$ points $(u_0,x_1(u_0)),\dots, (u_0,x_{d-2}(u_0))$ are very close to the point $(0,0)$, since the closure of the line $\{u=0\}$ meets the closure of the curve $E$ at infinity with multiplicity two and at the point $(0,0)$ with multiplicity $d-2$. Therefore
$$(\cup_{l=1}^d\mathcal U_{1,l})\cap h_d^{-1}(u_0)=\{ (u_0,x_s(u_0),y_l(u_0,x_s(u_0)))\}_{1\leq s\leq d-2, 1\leq l\leq d},$$
where $y_l(u_0,x_s(u_0))$, $l=1,\dots, d$, are the roots of the equation $$x^d_s(u_0)+ y^d +1+u_0x_s^2(u_0)=0,$$ and hence the intersection
$\mathcal U_{1,1}\cap h_d^{-1}(u_0)$ consists of $d-2$ different points for each $u_0\in L_{\varepsilon}\setminus \{ 0\}$. Note also that
$\mathcal U_{1,1}\cap h_d^{-1}(0)$ is the single point $(0,\overline z_{1,1})$. Therefore to prove that $\mathcal U_{1,1}$ is irreducible, it suffices to show that $\mathcal U_{1,1}\cap h_d^{-1}(L_{\varepsilon})$ is a smooth curve, since overwise the curve
$\mathcal U_{1,1}\cap h_d^{-1}(L_{\varepsilon})$ is the union of several components lying in different irreducible components of
$\mathcal U_{1,1}$ and meeting at the point $(0,\overline z_{1,1})$ which must be the singular point of
$\mathcal U_{1,1}\cap h_d^{-1}(L_{\varepsilon})$.

To prove the smoothness of $\mathcal U_{1,1}\cap h_d^{-1}(L_{\varepsilon})$  at $(0,\overline z_{1,1})$ note that in non-homogeneous coordinates $(u,x,y_1=y-\mu_1)$ the curve $\mathcal U_{1,1}\cap h_d^{-1}(L_{\varepsilon})$ is the complete intersection of two surfaces given by  equation (\ref{curve4}) and the equation $x^d+(y_1+\mu_1)^d+1+ux^2=0$. It is easy to check that these two surfaces are non-singular at the point $(u,x,y_1)=(0,0,0)$ and meet transversally at this point. \qed

\begin{cor} \label{cor1} Let $\mathcal U_{j,l}\subset \mathcal U_{\varepsilon}$ be the same as in Claim \ref{cl1}. Then $\mathcal U_{j,l}\setminus h_d^{-1}(\mathcal B_d)$ is a connected smooth variety.
\end{cor}

\subsection{Transitivity of the actions of the groups $\mathcal G_d$.} \label{trans}
In notation used in subsection \ref{hF}, let the base point $\overline a_0$ of the fundamental group
$\pi_1(\mathbb P^{K_d}\setminus \mathcal B_d,\overline a_0)$ lie in the neighbourhood $U_{\varepsilon}$, $\varepsilon <<1$. Then the set $I_{\overline a_0}=h_d^{-1}(\overline a)$ naturally splits into the union of $3d$ subsets,
$I_{\overline a_0}=\bigsqcup_{j=1}^3\bigsqcup_{l=1}^d I_{j,l}$, where
$I_{j,l}=I_{\overline a}\cap \mathcal U_{j,l}=\{ p_{j,l,1},\dots p_{j,l,d-2}\}$.

\begin{claim} \label{cl3} The group $\mathcal G_d=Im\, h_{d*}\subset \mathbb S_{3d(d-2)}$ acts transitively on the set $I_{\overline a_0}$.
\end{claim}

\proof In the beginning, let us show that the group $\mathcal G_d$ acts transitively on each subset $I_{j,l}$. Indeed, by Corollary \ref{cor1}, for each pair $(j,l)$ the variety $\mathcal U_{j,l}\setminus h_d^{-1}(\mathcal B_d)$ is connected. Therefore
any two points $p_{j,l,m_1}, p_{j,l,m_2}\in I_{j,l}$ can be connected by a smooth path $\gamma\subset \mathcal U_{j,l}\setminus h_d^{-1}(\mathcal B_d)$. Then the image $g=h_{d*}(\overline \gamma)\in \mathbb S_{3d(d-2)}$ of the element $\overline \gamma\in \pi_1(\mathbb P^{N_d}\setminus \mathcal B_d,\overline a)$ represented by the loop $h_d(\gamma)$ sends the point $p_{j,l,m_1}$ to $p_{j,l,m_2}$.

Now to complete the proof of Claim \ref{cl3}, it suffices to show that for each pair $(j,l)$ the point $p_{1,1,1}\in I_{1,1}$ can be connected with some point $p_{j,l,m}\in I_{j,l}$ by a smooth path $l\subset \mathcal I_d\setminus \mathcal B_d$. For this let us consider an element $g_1\in G_d\subset PGL(3,\mathbb C)$ such that $g_1(\mathcal U_{1,1})=\mathcal U_{j,l}$, where the group $G_d$ was introduced in Subsection \ref{hF}. The group $PGL(3,\mathbb C)$ is connected. Therefore we can find a smooth path $g_t\subset PGL(3,\mathbb C)$, $t\in [0,1]$, connecting the elements $g_0=Id$ and $g_1$ in $PGL(3,\mathbb C)$. Then the path $g_t(p_{1,1,1})$ lies in  $\mathcal I_d\setminus h_d^{-1}(\mathcal B_d)$, since for each $t\in [0,1]$ the curve $g_t(C_{\overline a})$ is smooth and it has not multiple inflection points, and the point $g_t(p_{1,1,1})$ is an inflection point of the curve $g_t(C_{\overline a})$. Note also that the path $g_t(p_{1,1,1})$ connects the point $p_{1,1,1}$ with some point $g_1(p_{1,1,1})\in \mathcal U_{j,l}$. As above, by Corollary \ref{cor1}, the point $g_1(p_{1,1,1})$ can be connected with any point lying in $I_{j,l}$ by a path in $\mathcal U_{j,l}\setminus h_d^{-1}(\mathcal B_d)$. \qed

\subsection{Behaviour of the covering $h_d$ near a $2$-tuple inflection point.} \label{hF2}
Let $p$ be a $2$-tuple inflection point of a curve $C$ of degree $d\geq 4$. Without loss of generality, we can assume that $p=(0,0,1)$ and the line $\{ z_1=0\}$ is the tangent line to $C$ at the point $p$. Then $C$ is given by equation $F(z_1,z_2,z_3)=0$, where $F(z_1,z_2,z_3)$ is a homogeneous polynomial of the form $z_2^4R(z_2,z_3)+z_1S(z_1,z_2,z_3)$ and where $R(z_1,z_3)$ is a homogeneous polynomial of degree $d-4$ such that $R(0,1)=1$, and $S(z_1,z_2,z_3)$ is a homogeneous polynomial of degree $d-1$ such that $S(0,0,1)=1$.

Let $V\subset \mathbb \mathbb P^{K_d}$ be a very small neighbourhood of the point $c$ corresponding to the curve $C$ and
$\mathcal V=h_d^{-1}(V)\subset \mathbb P^{K_d}\times \mathbb P^2$. Then $\mathcal V$ is the disjoint union of two components,
$\mathcal V=\mathcal V_1\bigsqcup \mathcal V_2$, where $\mathcal V_1$ contains the point $(c,p)\in \mathbb P^{K_d}\times \mathbb P^2$. Obviously, the restriction of $h_d$ to $\mathcal V_1$ has degree two, since $p$ is a $2$-tuple inflection point of $C$ and therefore under small deformation of $C$ the deformed curves,  near the point $p$, have either two different inflection points or one $2$-tuple inflection point.

\begin{claim}\label{2-tuple} There is a smooth curve $E_1\subset \mathcal V_1$ passing through the point $(c,p)$ and such that the restriction $h_{d|E_1}:E_1\to L_1$ of $h_d$ to $E_1$ is ramified at $(c,p)$ with multiplicity $\deg h_{d|E_1}=2$.
\end{claim}

\proof Consider the family of curves $C_v\subset \Delta\times \mathbb P^2\subset\mathbb P^{K_d}\times \mathbb P^2$  given by $F(z_1,z_2,z_3)+vz_2^2z_3^{d-2}=0$, where $\Delta =\{ |v|< \varepsilon_1 \}$ is the disk in $\mathbb C$ of small radius $\varepsilon_1$.

Consider the curve $E_1=\mathcal V_1\cap C_v\cap H_v$, were $H_v$ is the Hessian of the family $C_v$. Let  $x=z_1/z_3$ and $y=z_2/z_3$. Denote by $R'_j=\frac{\partial R}{\partial z_j}(x,1)$,
$S'_j=\frac{\partial S}{\partial z_j}(x,y,1)$, $R''_{j,l}=\frac{\partial^2 R}{\partial z_j\partial z_l}(x,1)$, and $R''_{j,l}=\frac{\partial^2 R}{\partial z_j\partial z_l}(x,y,1)$. In the coordinates $(v,x,y)$ the family $C_v$ is given by
\begin{equation} \label{a} y^4R(y,1)+vy^2+xS(x,y,1)=0\end{equation}
and the Hessian $H_v$ is given by
\begin{equation} \label {b}  { \tiny
\left|
 \begin{matrix} 2S'_1+xS''_{1,1}, & S'_2+xS''_{1,2}, & S'_3+xS''_{1,3}
\\
S'_2+xS''_{1,2}, & 12y^2R+8y^3R'_2+y^4R''_{2,2}+2v+xS''_{2,2}, & 4y^3R'_3+y^4R''_{2,3}+2\alpha vy+xS''_{2,3} \\
S'_{3}+xS''_{1,3}, & 4y^3R'_3+y^4R''_{2,3}+2\alpha vy+xS''_{2,3} & y^4R''_{3,3}+\beta vy^2 +xS''_{3,3}
\end{matrix} \right|=  0,}
\end{equation}
where $\alpha =d-2$ and $\beta=(d-2)(d-3)$.

Let $S'_j(0,0,1)=s_j$ and $S''_{j,l}(0,0,1)=s_{j,l}$. By assumption, $S(0,0,1)=1$. Therefore $S'_3(0,0,1)=d-1$ and it is easy to see that the differentials at the point $(v,x,y)=(0,0,0)$ of the polynomials in the left side of equations (\ref{a}) and (\ref{b}) are equal respectively to $\overline dx$ and $(2s_2s_3s_{2,3}-s_2^2s_{3,3})\overline dx-2(d-1)^2\overline dv$ (here we denote by $\overline d f$ the differential of function $f$ in order not to confuse with degree $d$ of the curve $C$). Therefore the surfaces $C_v$ and $H_v$ are nonsingular at the point $(0,0,0)$ and meet transversally at this point, since these differentials are linear independent. It follows from this that
$E_1=\mathcal V_1\cap C_v\cap H_v$ is a smooth curve and the differential of $h_{d|E_1}:E_1 \to \Delta$ vanishes at the point $(0,0,0)$. Therefore $\deg h_{d|E_1}=2$ since if a holomorphic surjective map of a smooth curve has degree one, then its differential vanishes nowhere. \qed

\begin{cor} \label{cor2} Let $C$ be a smooth curve of degree $d\geq 4$ having $3d(d-2)-1$ inflection points. Then, in notations used in the proof of Claim {\rm \ref{2-tuple}}, $h_d^{-1}(\Delta)\cap \mathcal V$ is the disjoint union of $3d(d-2)-1$ smooth curves, $h_d^{-1}(\Delta)\cap \mathcal V=\bigsqcup_{j=1}^{3d(d-2)-1}E_j$, where $h_{d|E_j}:E_j\to \Delta$ is a bi-holomorphic map for $j=2,\dots, 3d(d-2)-1$ and $h_{d|E_{1}}$ is a two-sheeted covering branched at the point $\{ v=0\}\in L_1$.
\end{cor}
Consider a point $(\overline a,\overline z)\in h_d^{-1}(\mathfrak M_d)$, where $\mathfrak M_d\subset\mathbb P^{K_d}$ is the variety defined in Subsection \ref{beg}. Lemma \ref{non-si} and Claim \ref{2-tuple} imply

\begin{prop} Let $\overline a_0\in \mathfrak M_d$ and $\overline z_0$ be a $2$-tuple inflection point of the curve $C_{\overline a_0}$.
Then $(\overline a_0,\overline z_0)$ is a smooth point of the variety $\mathcal I_d$.
\end{prop}

Corollary \ref{cor2} and Lemma \ref{pi} imply
\begin{claim} \label{transposition} For $d\geq 4$ the group $\mathcal G_d\subset \mathbb S_{3d(d-2)}$ contains a transposition.
\end{claim}

\subsection{Case $d=3$}
Consider the {\it Hesse pencil}, that is, the one-dimensional linear system of plane cubic curves given  by
\begin{equation}\label{Fer3} C_{(t_1,t_2)}: \quad t_1(z_1^3+z_2^3+z_3^3)+ t_2z_1z_2z_3=0, \quad (t_1,t_2)\in \mathbb P^1 ,\end{equation}
We call the surface  $\mathcal H\subset L\times \mathbb P^2\subset \mathbb P^{K_3}\times \mathbb P^2$ given in $L\times \mathbb P^2$ by equation (\ref{Fer3}) the {\it body} of the Hesse pencil, where $L\simeq \mathbb P^1$ and $K_3=9$.

It is easy to see that $\mathcal H$ is a smooth surface and the restriction $\sigma :\mathcal H\to\mathbb P^2$ of $\text{pr}_2$ to
$\mathcal H$  is the composition of nine $\sigma$-processes of $\mathbb P^2$ with centers at the base points of the Hesse pencil. Let $E_{q_j}=\sigma^{-1}(q_j)$, $j=1,\dots,9$, be the exceptional curve of $\sigma$ over the base point $q_j\in \mathbb P^2$ of the pencil.
The curves $E_j$ are sections of the projection $\text{pr}_1:L\times\mathbb P^2\to L$.

It is well-known (see, for example, \cite{BK}) that the base points of the Hesse pencil are the inflection points of each smooth member of the pencil.
The Hesse pencil has four degenerate members, $C_{(0,1)}$, $C_{(1,-3)}$, $C_{(1,-3e^{2\pi i/3})}$, and $C_{(1,-3e^{4\pi i/3})}$. Each of the degenerate members is the union of three lines.
Therefore
$$\mathcal I_3\cap \mathcal H=C_{(0,1)}\cup C_{(1,-3)}\cup C_{(1,-3e^{2\pi i/3})}\cup C_{(1,-3e^{4\pi i/3})}\cup(\cup_{j=1}^9E_j).$$

The group ${Hes}\subset PGL(3,\mathbb C)$ of the projective transformations leaving invariant the set of the inflection points of the Fermat curve $F_3=C_{(1,0)}$ is well investigated (see, for example, \cite{BK}). The order of ${Hes}$ is equal to $216$ and the action of ${Hes}$ on the $9$ inflection points of $F_3$ defines an imbedding ${Hes}\subset \mathbb S_9$ such that ${Hes}$ is a $2$-transitive  subgroup of
$\mathbb S_9$. The orbit of the Fermat curve $F_3$ under the action of ${Hes}$ consists of four members of the Hesse pencil: $F=C_{(1,0)},
C_{(1,6)}, C_{(1,6e^{2\pi i/3})}, C_{(1,6e^{4\pi i/3})}$. We choose three continuous paths $l_j$, $j=0,1,2$, in $L\setminus \{ (0,1),(1,-3),(1,-3e^{2\pi i/3}),(1,-3e^{4\pi i/3})\}$ connecting the points $(1,6e^{2j\pi i/3})$ with the point $(1,0)$.

It is well-known (\cite{He1}, \cite{He2}) that the set of nine inflection points of a plane cubic curve is a projectively
rigid set, that is, for each two smooth plane cubic curves $C_1$ and $C_2$ there is a projective transformation of the plane sending the set of the inflection points of $C_1$ onto the  set of the inflection points of $C_2$. Therefore there is an imbedding $\varphi : PGL(3,\mathbb C)\to (\mathbb P^2)^9$  given for $\tau\in PGL(3,\mathbb C)$ by
$$\varphi: \tau\mapsto (\tau(q_1),\dots,\tau(q_9))\in (\mathbb P^2)^9,$$
where $\{ q_1,\dots,q_9\}\subset \mathbb P^2$ is the set of the inflection points of the Fermat curve $F_3$. (Note that $\varphi$ depends on the numbering of the inflection points of the Fermat curve $F_3$, that is, there are $\frac{9!}{216}$ such imbeddings.) Denote by $\mathcal P=\varphi (PGL(3,\mathbb C))\subset (\mathbb P^2)^9$.

\begin{claim} \label{case3} We have $\mathcal G_3=Hes\subset \mathbb S_9$.
\end{claim}
\proof
Consider the homomorphism $h_{3*}:\pi_1(\mathbb P^{K_3}\setminus \mathcal S_3,f)\to \mathbb S_9$, where $f =(1,0)\in L\subset \mathbb P^{K_3}$ is the point corresponding to the Fermat curve $F_3$. Then $\mathcal G_3=h_{3*}(\pi_1(\mathbb P^{K_3}\setminus \mathcal S_3,f))$ acts on the set $h_{3}^{-1}(f)=\{ q_1, \dots, q_9\}$.

Let us show that $Hes\subseteq \mathcal G_3$. For this, consider an element $g_1\in Hes\subset PGL(3,\mathbb C)$. Since $PGL(3,\mathbb C)$ is a connected variety, we can choose a  continuous path $g_t\subset PGL(3,\mathbb C)$, $0\leq t\leq 1$, connecting $g_0=Id$ with $g_1$.
The path $g_t$ defines a continuous path $l(t)\subset \mathbb P^{K_3}\setminus \mathcal S_3$ such that $C_{l(t)}=g_t(F_3)$. We have $C_{l(1)}$ is a member of the Hesse pencil, since $g_1\in Hes$. Denote by $\Gamma\subset \mathbb P^{K_3}\setminus \mathcal S_3$ the path $l(t)$ if $C_{l(1)}=F_3$ and $l(t)\cup l_j$ if $C_{l(1)}=C_{(1,6e^{2j\pi i/3})}$. Then the loop $\Gamma$ represents an element $\gamma \in\pi_1(\mathbb P^{K_3}\setminus \mathcal S_3,f)$ and it is easy to see that the action of $h_{3*}(\gamma)$ on $h_{3}^{-1}(f)=
\{ q_1,\dots,q_9\}$ is the same as the action of $g_1$.

Let us show that $\mathcal G_3\subseteq Hes$. For this, consider an element $g\in\mathcal G_3$ and  a continuous loop $\Gamma(t)\subset \mathbb P^{K_3}\setminus \mathcal S_3$ starting and ending at $f$ and representing an element $\gamma \in \pi_1(\mathbb P^{K_3}\setminus \mathcal S_3,f))$ such that $h_{3*}(\gamma)=g$. We lift the loop $\Gamma(t)$  to $\mathcal I_3$ and this lift consists of $9$ continuous  paths $\Gamma_1(t),\dots,\Gamma_9(t)\subset \mathcal I_3$  starting and ending at the points of $h_{3*}^{-1}(f)$.  So, we obtain $9$ continuous paths $\text{pr}_2(\Gamma_1(t)),\dots,\text{pr}_2(\Gamma_9(t))\subset \mathbb P^2$. If we number the paths $\Gamma_j(t)$ so that $\Gamma_j(0)=q_j$, then $(\text{pr}_2(\Gamma_1(t)),\dots,\text{pr}_2(\Gamma_9(t)))$ is a continuous path in $\mathcal P$, since $\{\text{pr}_2(\Gamma_1(t)),\dots,\text{pr}_2(\Gamma_9(t))\}$ is the set of the inflection points of smooth plane cubic curves. Therefore there is an element $\tau\in PGL(3,\mathbb C)$ such that $\tau (q_j)=\text{pr}_2(\Gamma_j(1))\in h_3^{1}(f)$ for $j=1,\dots, 9$, that is, $g=\tau\in Hes$. \qed

\subsection{Case $d=4$.}
Consider the Klein curve $Kl\subset\mathbb P^2$ given by
$$z_1^3z_2+z_1^3z_3+z_1z_2^3+z_2^3z_3+z_1z_3^3+z_2z_3^3=0.$$
It is well-known (see, for example, \cite{D}) that the automorphism group $Aut(Kl)$ of $Kl$ have the following properties.
The order of $Aut(Kl)$ is equal to $168$ and $Aut(Kl)\simeq PSL(2,\mathbb Z_7)$, the group $Aut(Kl)$ can be represented as a subgroup of
$PGL(3,\mathbb C)$ leaving invariant the curve $Kl$ and the set of inflection points is an orbit under the action $Aut(Kl)$, the order of the stabilizer of each inflection point is equal to $7$. In particular, the action of $Aut(Kl)$ on the set of the inflection points of $Kl$ is transitive.
\begin{claim} \label{K} There is an imbedding $Aut(Kl)\subset \mathcal G_4\subset \mathbb S_{24}$ such that $Aut(Kl)$ is a transitive subgroup of $\mathbb S_{24}$.
\end{claim}
\proof Let $\overline a_0\in \mathbb P^{K_4}$ be the point corresponding to the curve $Kl$ and $g_1$ an element of $Aut(Kl)\subset PGL(3,\mathbb C)$.
Since $PGL(3,\mathbb C)$ is connected, there is a continuous path $g_t\subset PGL(3,\mathbb C)$, $t\in [0,1]$, connecting $g_0=Id$ and $g_1$. Then the loop $\Gamma\subset \mathbb P^{K_4}\setminus \mathcal B_4$ given by $g_t(\overline a_0)$, $t\in [0,1]$, defines  an element
$\gamma\in \pi_1(\mathbb P^{K_4}\setminus \mathcal B_4,k)$ such that the action of
$h_{4*}(\gamma)\in \mathbb S_{24}$ on the set $h_4^{-1}(\overline a_0)$ of the inflection points of $Kl$ coincides with the action of $g_1\in Aut(Kl)$. \qed \\

By Claim  \ref{imbedding} and  \ref{K}, the group $\mathcal G_4\subset\mathbb S_{24}$ has the following properties:
\begin{itemize}
\item[(1)] $\mathcal G_4$ contains a subgroup $Aut (Kl)$ which acts transitively on the set $I_{24}=\{ 1,2,\dots,23,24\}$;
\item[(2)] there are a partition $I_{24}=J_1\bigsqcup J_2$, $|J_1|=9$, $|J_2|=15$, and a quasi-imbedding $\mathcal G_3\prec\mathcal G_4$ such that $J_1$ is invariant under the action of
    $\widetilde{\mathcal G_3}\subset \mathcal G_4$ {\rm (see subsection \ref{quasi})} and the action of $\widetilde{\mathcal G_3}$ on $J_1$ is $2$-transitive;
\item[(3)] the group $\mathcal G_4$ contains a transposition.
\end{itemize}

\begin{claim} \label{g4}
Properties $(1)$ -- $(3)$ imply $\mathcal G_4=\mathbb S_{24}$.
\end{claim}
\proof We say that a subset $J\subset I_{24}$ is $2$-{\it transitive} (with respect to the action of $\mathcal G_4$) if for each two pairs $\{ j_1,j_2\}\subset J$ and $\{ j_3,j_4\}\subset J$ there is an element $g\in \mathcal G_4$ such that $g(\{ j_1,j_2\})=\{ j_3,j_4\}$.

Denote by $\widetilde J\subset I_{24}$ a $2$-transitive subset of maximal cardinality. Obviously, if $J\subset I_{24}$ is a $2$-transitive subset then for each $g\in \mathcal G_4$ the set $g(J)$ is also $2$-transitive, and if $J_1$ and $J_2$ are  $2$-transitive subsets such that the cardinality $|J_1\cap J_2| \geq 2$, then $J_1\cup J_2$ is also $2$-transitive. Therefore it is easy to see that there are two possibilities: either $\widetilde J=I_{24}$, or $|\widetilde J|=12$, since, by property (2), the cardinality $|\widetilde J|\geq 9$ and, by property (1), $\mathcal G_4$ acts transitively on $I_{24}$.

Let us show that the second case is impossible. Indeed, in this case it follows from transitivity of the action of $\mathcal G_4$ that for each $g\in\mathcal G_4$ either $g(\widetilde J)=\widetilde J$, or $g(\widetilde J)=I_{24}\setminus \widetilde J$. Therefore the action of $\mathcal G_4$ on $I_{24}$ induces an action on the set $\{ \widetilde J, I_{24}\setminus \widetilde J \}$ of cardinality two,  that is, there is an epimorphism $\varphi :\mathcal G_4\to \mathbb Z_2$. But, in this case the restriction  $\varphi_{|Aut(Kl)}: Aut(Kl)\to \mathbb Z_2$ is also an epimorphism, since $Aut(Kl)$ acts transitively on the set $I_{24}$. On the other hand, $Aut(Kl)\simeq PSL(2,\mathbb Z_7)$ is a simple group and therefore $\varphi_{|Aut(Kl)}$ must be trivial homomorphism.

Now, to complete the proof of Claim \ref{g4}, it suffices to apply property (3), since $\mathcal G_4$ acts $2$-transitively on $I_{24}$ and therefore the all transpositions of $\mathbb S_{24}$ are contained in $\mathcal G_4$. \qed

\subsection{The end of the proof of Theorem \ref{main}}
To complete the proof of Theorem \ref{main} we need in the following
\begin{lem} \label{act} Let G be a subgroup of the symmetric group $\mathbb S_m$ acting on a finite set $M$ of cardinality $m$.
Assume that
\begin{itemize}
\item[$(i)$] $G$ acts transitively on $M$;
\item[$(ii)$] there are a subgroup $G_1$ of $G$ and a subset $M_1$ of $M$ 
such that
\begin{itemize}
\item[$(ii)_1$] $M_1$ is invariant under the action of the group $G_1$,
\item[$(ii)_2$] $2m_1\geq m+2$, where $m_1$ is the cardinality of the set $M_1$,
\item[$(ii)_3$] the action of the group $G_1$ on $M_1$ is $2$-transitive;
\end{itemize}
\item[$(iii)$] the group $G$ contains a transposition.
\end{itemize}
Then $G=\mathbb S_m$.
\end{lem}

\proof By $(ii)_3$, for each element $g\in G$ the subgroup $gG_1g^{-1}$ of $G$ acts $2$-transitively on $g(M_1)$ and by $(ii)_2$,
the group $G$ acts $2$-transitively on $M_1\cup g(M)$, since there are at least two elements in the intersection of $M_1\cap g(M_1)$.
It follows from $(i)$ that for each element $x\in M$ there is an element $g\in G$ such that $x\in g(M_1)$. Therefore $G$ acts $2$-transitively on $M$ and hence applying $(iii)$ the group $G\subset \mathbb S_m$ contains all transpositions. \qed

Now, applying induction on $d$, Theorem \ref{main} follows from Claims  \ref{imbedding}, \ref{cl3}, \ref{transposition}, \ref{case3},  \ref{g4} and Lemma \ref{act}, since
$$2[3d(d-2)]\geq [3(d+1)(d-1)]+2$$
for $d\geq 4$.

\section{Behaviour of the covering $h_d$ near  a node of a nodal curve} \label{sing}
\subsection{On the subset of $\mathcal S_d$ consisting of the points corresponding to the nodal curves.}
Denote by $\mathcal N_d$ a subvariety of $\mathcal S_d$ consisting of the points $\overline a\in \mathcal S_d$ such that the set of singular points of the curves $C_{\overline a}$ consist of the only one ordinary node. The following Proposition is well-known.

\begin{prop} \label{propx} The variety $\mathcal S_d$ is an irreducible hypersurface in $\mathbb P^{K_d}$ for each $d\geq 3$. The variety $\mathcal N_d$  is a non-empty Zariski open subset of $\mathcal S_d$.
\end{prop}
{\it Proof} of this proposition is similar to the proof of Proposition \ref{prop1} and therefore it will be omitted. \qed 

\begin{claim} \label{singS} The variety $\mathcal N_d\subset \mathbb P^{K_d}$ is smooth.
\end{claim}
\proof Consider a point $\overline a_0\in \mathcal N_d$.
Without loss of generality, we can assume that $\overline z_0=(0,0,1)$ is the singular point of $C_{\overline a_0}$ and $C_{\overline a_0}$ is given by equation
$F(\overline a_0,\overline z)=0$, where
$$F(\overline a_0,\overline z)=z_1z_2z_3^{d-2}+ R(z_1,z_2,z_3),$$
and where $R(z_1,z_2,z_3)$ is a homogeneous  polynomial of degree $d$ such that it, as a polynomial in variable $z_3$, has degree  $\leq d-3$. In particular, the coordinates $\alpha_{0,0,d}$, $\alpha_{1,0,d-1}$, $\alpha_{0,1,d-1}$, $\alpha_{2,0,d-2}$, $\alpha_{0,2,d-2}$ in $\overline a_0$ are equal to zero. Let us show that the tangent space $T_{\overline a_0}\mathcal S_d$ to $\mathcal S_d$ at $\overline a_0$ is given by equation $a_{0,0,d}=0$.

The variety $\mathcal C_d\subset \mathbb P^{K_d}\times \mathbb P^2$ at the point $(\overline a_0,\overline z_0)$ in non-homogeneous coordinates is given by
\begin{equation}\label{vars1} a_{0,0,d}+ a_{1,0,d-1}x+ a_{0,1,d-1}y +xy+a_{2,0,d-2}x^2+a_{0,2,d-2}y^2+
\widetilde R(x,y,1)=0. \end{equation}
It is easy to see that $\mathcal C_d$ is non-singular at $(\overline a_0,\overline z_0)$ and the elements of the set of variables
$\{ x,y\}\cup \{a_{k,m,n}\}_{k+m+n=d}\setminus \{ a_{0,0,d}, a_{1,1,d-2}\}$ are local coordinates at $(\overline a_0,\overline z_0)$.
Consider in $\mathcal C_d$ a subvariety $Sing_d$ given by
\begin{equation}\label{vars2}
a_{1,0,d-1}+y+2a_{2,0,d-2}x+\widetilde R'_x(x,y,1)=0
\end{equation}
and
\begin{equation}\label{vars3}
a_{0,1,d-1}+x+2a_{0,2,d-2}y+\widetilde R'_y(x,y,1)=0.
\end{equation}
Obviously, in a neighbourhood of $\overline a_0$ the variety $\mathcal S_d$ is the image of $Sing_d$ under the morphism $\widetilde h_d$. It follows from (\ref{vars1}) -- (\ref{vars3}) that $Sing_d$ is non-singular at
$(\overline a_0,\overline z_0)$, the elements of the set of variables $\{a_{k,m,n}\}_{k+m+n=d}\setminus \{ a_{0,0,d}, a_{1,1,d-2}\}$  are local coordinates at $(\overline a_0,\overline z_0)$ and $\widetilde h_d: Sing_d\to\mathcal S_d$ is given in these coordinates by
$$a_{0,0,d}= a_{1,0,d-1}^2A_1(a_{k,m,n})+a_{1,0,d-1}a_{0,1,d-1}A_2(a_{k,m,n})+a_{0,1,d-1}^2A_3(a_{k,m,n}),$$
where $A_j(a_{k,m,n})$, $j=1,2,3$, are power serieses in variables
$$\{a_{k,m,n}\}_{k+m+n=d}\setminus \{ a_{0,0,d}, a_{1,1,d-2}\}. \qquad \qed $$

Let $\overline z_0$ be an ordinary node of $C_{\overline a_0}$. In \cite{Ku}, it was shown that $(C,H_C)_{\overline z_0}=6$ if $\overline z_0$ is not an inflection point of each of the branches of $C_{\overline a_0}$ at $\overline z_0$. Therefore the local degree of $h_d$ at the point $(\overline a,\overline z_0)$ is equal to $\deg_{(\overline a_0,\overline z_0)} h_d=6$.

\subsection{On the local monodromy groups of $h_d$ at the points corresponding to the nodal curves.} Denote $\mathfrak N_d=\mathcal N_d\setminus \mathcal M_d$.

\begin{prop} \label{propx} Let $\overline a_0\in \mathfrak N_d$ and $\overline z_0$ be the singular point of $C_{\overline z_0}$. Then the local monodromy group $\mathcal G_{d,\overline a_0}\subset \mathbb S_{3d(d-2)}$ at the point $\overline a_0$ is a cyclic group of order $3$ and it is generated by the product of two disjoint cycles of length $3$.
\end{prop}
\proof
Without loss of generality we can asume that  $C_{\overline a_0}$ is given by equation
$F(\overline a_0,\overline z)=0$, where
$$F(\overline a_0,\overline z)=z_1z_2z_3^{d-2}+ z_3^{d-3}\sum_{j+k=3}\alpha_{j,k,d-3}z_1^jz_2^k+
R(z_1,z_2,z_3)$$
and where $R$ is a  polynomial of degree  $\leq d-4$ in the variable $z_3$.

Note that $a_{3,0,d-3}\neq 0$ and $a_{0,3,d-3}\neq 0$ if $\overline a_0\in \mathfrak N_d$. Indeed, if, for example, $a_{3,0,d-3}= 0$ then
it is easy to see that the line $L$ given by $t\overline u+\overline a_0$, where $t\in \mathbb C$ and in $\overline u$ all coordinates except the coordinate $u_{0,1,d-1}$ are equal to zero and $u_{0,1,d-1}=1$, lies in $\mathcal M_d$.

Consider a one-parametric family $C_{\overline a_t}$ given by equation
$$F(\overline a_0,\overline z)+tz_3^d=0$$
and its projection $\text{pr}_1(C_{\overline a_t}) =L=\{ \overline a_t=\overline a_0+ t\overline v\}\subset \mathbb P^{K_d}$, where in $\overline v$ all coordinates except the coordinate $v_{0,0,d}$ are equal to zero and $v_{0,0,d}=1$. By Claim \ref{singS}, $L$ meets $\mathcal S_d$ transversally at $\overline a_0$.

In non-homogeneous coordinates $x=\frac{z_1}{z_3}, y=\frac{z_2}{z_3}$ we have $\overline z_0=(0,0)$, the family $C_{\overline a_t}$ in
$L\times\mathbb C^2\subset\mathbb P^{K_d}\times\mathbb P^2$ is given by equation
\begin{equation} \label{node1} t+xy+ \sum_{i+j=3}a_{i,j,d-3}x^iy^j+ \text {terms of higher degree}=0,\end{equation}
and everybody can easily check that its Hessian $H_{C_{\overline a_t}}$  is given by equation
\begin{equation} \label{node2}
\begin{array}{l} 2(d-2)^2(xy-3a_{3,0,d-3}x^3+a_{2,1,d-3}x^2y+a_{1,2,d-3}xy^2-3a_{0,3,d-3}y^3)+ \\ d(d-1)(1+4a_{2,1,d-3}x+4a_{1,2,d-3}y)t+r_1(x,y)+tr_2(x,y)=0,
\end{array}\end{equation}
were $r_1(x,y)=\sum_{i+j\geq 4}b_{i,j}x^iy^j$ and $r_2(x,y)=\sum_{i+j\geq 2}c_{i,j}x^iy^j$ are some polynomials.

Consider the curve $Z=h_d^{-1}(L)=C_{\overline a_t}\cap H_{C_{\overline a_t}}\subset X$, where $X$ is a surface in $\mathbb C^3$ given by equation (\ref{node1}). It is easy to see that $X$ is isomorphic to $\mathbb C^2=Spec\, \mathbb C[x,y]$ and $Z$ in $X$ is given by equation
\begin{equation} \label{node3}
\begin{array}{l} (d^2-7d+8)xy-6(d-2)^2(a_{3,0,d-3}x^3+a_{0,3,d-3}y^3) - \\
2(d^2+2d-4)(a_{2,1,d-3}x^2y+a_{1,2,d-3}xy^2)+\text {terms of higher degree}=0, \end{array}
\end{equation}
since
\begin{equation}\label{t}  t=-(xy+ \sum_{i+j=3}a_{i,j,d-3}x^iy^j+ \text {terms of higher degree}).\end{equation}

It follows from (\ref{node3}) that $Z$ has a node at the point ${\bf p}=(\overline a_0,\overline z_0)$. To resolve this point, consider the $\sigma$-process $\sigma: \widetilde X\to X$ with center at ${\bf p}$. The surface $\widetilde X$ is covered by two open neighbourhoods isomorphic to $\mathbb C^2$, $\widetilde X=U_1\cup U_2$. The coordinates in $U_j$, $j=1,2$, are $x_j$, $y_j$ and $\sigma_{\mid U_1}$ is given by $x=x_1$ and $y=x_1y_1$, and $\sigma_{\mid U_2}$ is given by $x=x_2y_2$ and $y=y_2$. Therefore the proper inverse image
$\sigma^{-1}(Z)\cap U_1$ is given by equation
\begin{equation} \label{node5}
\begin{array}{l} (d^2-7d+8)y_1-6(d-2)^2(a_{3,0,d-3}x_1+a_{0,3,d-3}x_1y_1^3) - \\
2(d^2+2d-4)(a_{2,1,d-3}x_1y_1+a_{1,2,d-3}x_1y_1^2)+\text {terms of higher degree}=0, \end{array}
\end{equation}
and therefore the curve $\widetilde Z=\sigma^{-1}(Z)$ is non-singular at the point ${\bf p}_1=\widetilde Z\cap U_1\cap E$, where
$E$ is the exceptional divisor of $\sigma$.

Since $a_{3,0,d-3}\neq 0$, the coordinate $x_1$ is a local parameter in $\widetilde L$ at the point ${\bf p}_1$ and $$y_1=\frac{6(d-2)^2a_{3,0,d-3}}{d^2-7d+8}x_1+\sum_{j=2}^{\infty} b_jx_1^j.$$

It follows from (\ref{t}) that
$$\begin{array}{lcll} t = & x_1^2y_1+ a_{3,0,d-3}x_1^3 & + & \text {terms of higher degree}= \\
& (\frac{6(d-2)^2}{d^2-7d+8}+1)a_{3,0,d-3}x_1^3 & + & \text {terms of higher degree.} \end{array}$$
Therefore the covering $h_d\circ\sigma:\widetilde Z\to L$ is ramified at ${\bf p}_1$ with multiplicity three.

Similar calculations (which will be omitted) show that the covering $h_d\circ\sigma:\widetilde Z\to L$ also is ramified at
${\bf p}_2=\widetilde Z\cap U_2\cap E$ with multiplicity three, since $a_{0,3,d-3}\neq 0$. \qed \\

Let $\nu_d:\mathfrak I_d\to \mathcal I_d$ be the normalization of the variety $\mathcal I_d$ and
$\overline h_{d}=h_d\circ\nu_d:\mathfrak I_d\to \mathbb P^{K_d}$. The following Proposition is an easy corollary of Lemma \ref{pi}, Claim \ref{locgroup}, and Proposition \ref{propx}.
\begin{prop} Let $\overline a_0\in \mathfrak N_d$ and $\overline z_0$ is the singular point of the curve $C_{\overline a_0}$. Then
\begin{itemize}
\item[$(i)$] the variety $\mathfrak I_d$ is smooth at the point ${\bf p}=(\overline a_0,\overline z_0)$;
\item[$(ii)$] $\nu_d^{-1}({\bf p})=\{ {\bf q}_1, {\bf q}_2\}$ consists of two points;
\item[$(iii)$] $\overline h_d$ is ramified along $\nu_d^{-1}(\mathfrak N_d)$ and the local degree $\deg_{{\bf q}_j}\overline h_d=3$ for $j=1,2$.
\end{itemize}
\end{prop}

\ifx\undefined\bysame
\newcommand{\bysame}{\leavevmode\hbox to3em{\hrulefill}\,}
\fi

\ifx\undefined\bysame
\newcommand{\bysame}{\leavevmode\hbox to3em{\hrulefill}\,}
\fi

\end{document}